\documentclass[10pt]{amsart}
\usepackage{amsfonts}
\usepackage{ifthen}
\usepackage{amsthm}
\usepackage{amsmath}
\usepackage{graphicx}
\usepackage{amscd,amssymb,amsthm}
\usepackage{graphicx}
\usepackage{epstopdf}
\usepackage{hyperref}

\newcounter{minutes}
\divide\time by 60
\newcounter{hours}
\multiply\time by 60 \addtocounter{minutes}{-\time}

\newtheorem{theorem}{Theorem}

\keywords{Convex functions; radius of convexity; Mittag-Leffler expansions; $q$-Bessel functions; zeros of $q$-Bessel functions; Laguerre-P\'olya class of entire functions.}
\subjclass[2010]{30C45, 30C15, 33C10}

\title{Bounds for radii of convexity of some $q$-Bessel functions}

\author[H. Orhan]{Hal\.{i}t Orhan}
\address{Department of Mathematics, Faculty of Science, Atat\"{u}rk University, Erzurum, Turkey}
\email{orhanhalit607@gmail.com}
\author[\.{I}. Akta\c{s}]{\.{I}brah\.{I}m Akta\c{s}}
\address{Department of Mathematical Engineering, Faculty of Engineering and Natural Sciences, G\"{u}m\"{u}\c{s}hane University, G\"{u}m\"{u}\c{s}hane, Turkey}
\email{aktasibrahim38@gmail.com}

\begin{document}

\def\thefootnote{}
\footnotetext{ \texttt{File:~\jobname .tex,
          printed: \number\year-\number\month-\number\day,
          \thehours.\ifnum\theminutes<10{0}\fi\theminutes}
} \makeatletter\def\thefootnote{\@arabic\c@footnote}\makeatother

\maketitle

\begin{abstract}
In the present investigation, by applying two different normalizations of the Jackson and Hahn-Exton $q$-Bessel functions tight lower and upper bounds for the radii of convexity of the same functions are obtained. In addition, it was shown that these radii obtained are solutions of some transcendental equations. The known Euler-Rayleigh inequalities are intensively used in the proof of main results. Also, the Laguerre-P\'olya class of real entire functions plays an important role in this work.
\end{abstract}

\section{Introduction}
Let $\mathbb{D}_r$ be the open disk $\{z\in\mathbb{C}:\left|z\right|<r\}$ with radius $r>0$ and $\mathbb{D}_1=\mathbb{D}$. Let $\mathcal{A}$ denote the class of analytic functions $f:\mathbb{D}_r\rightarrow\mathbb{C},$
$$f(z)=z+\sum_{n\geq2}a_{n}z^n,$$
which satisfy the normalization conditions $f(0)=f^{\prime}(0)-1=0$. By $\mathcal{S}$ we mean the class of functions belonging to $\mathcal{A}$ which are univalent in $\mathbb{D}_r$. The class of convex functions is defined by $$\mathcal{K}=\bigg\{f\in \mathcal{S}:\Re\left(1+\frac{zf^{\prime \prime}(z)}{f^{\prime}(z)}\right)>0 \text{ for all } z\in \mathbb{D}_r \bigg\}.$$ It is known that the convex functions do not need to be normalized, namely, the definition of $\mathcal{K}$ is also valid non-normalized analytic function $f:\mathbb{D}\rightarrow\mathbb{C}$ which has the property $f^{\prime}(0)\neq0$. The radius of convexity of an analytic locally univalent function $f:\mathbb{C}\rightarrow\mathbb{C}$ is defined by $$r^{c}(f)=sup\bigg\{r>0:\Re\left(1+\frac{zf^{\prime \prime}(z)}{f^{\prime}(z)}\right)>0 \text{ for all } z\in \mathbb{D}_r \bigg\}.$$ Note that $r^{c}(f)$ is the largest radius for which the image domain $f\left(\mathbb{D}_{r^{c}(f)}\right)$ is a convex domain in $\mathbb{C}.$ For more information about convex functions we refer to Duren's book \cite{Duren} and to the references therein.

The Jackson and Hahn-Exton $q$-Bessel functions are defined as follows:
$$J_{\nu}^{(2)}(z;q)=\frac{(q^{\nu+1};q)_{\infty}}{(q;q)_{\infty}}\sum_{n\geq0}\frac{(-1)^{n}\left(\frac{z}{2}\right)^{2n+\nu}}{(q;q)_{n}(q^{\nu+1};q)_{n}}q^{n(n+\nu)}$$ and $$J_{\nu}^{(3)}(z;q)=\frac{(q^{\nu+1};q)_{\infty}}{(q;q)_{\infty}}\sum_{n\geq0}\frac{(-1)^{n}z^{2n+\nu}}{(q;q)_{n}(q^{\nu+1};q)_{n}}q^{\frac{1}{2}{n(n+1)}},$$ where $z\in\mathbb{C},\nu>-1,q\in(0,1)$ and $$(a;q)_0=1,	(a;q)_n=\prod_{k=1}^{n}\left(1-aq^{k-1}\right),	 (a,q)_{\infty}=\prod_{k\geq1}\left(1-aq^{k-1}\right).$$
It is known that the Jackson and Hahn-Exton $q$-Bessel functions are $q$-extensions of the classical Bessel function of the first kind $J_{\nu}$. Clearly, for fixed $z$ we have $$J_{\nu}^{(2)}\left((1-z)q;q\right)\rightarrow{J}_{\nu}(z)$$ and $$J_{\nu}^{(3)}\left((1-z)q;q\right)\rightarrow{J}_{\nu}(2z)$$ as $q\nearrow1.$ The readers can find the properties of Jackson and Hahn-Exton $q$-Bessel functions in \cite{ismail1,ismail2,koelink,Koornwinder} and also comprehensive information on the Bessel function of the first kind can be found in Watson's treatise \cite{Wat}. Recently, the geometric properties of some special functions (like Bessel, Struve, Lommel and Wright functions) have been investigated by  many authors (see \cite{aktas2,aktas1,mathematica,publ,lecture,bdoy,bsk,BTK,bos,samy,basz,basz2,BY,BCDT,brown,ismail3,todd,wilf}). Also, the authors in \cite{aktas3,BDM} have studied the radii of starlikeness and convexity of some $q$-Bessel functions. In particular, tight lower and upper bounds for the radii of starlikeness of some $q$-Bessel functions were obtained in \cite{aktas3}. Most of above papers benefited from some properties of the positive zeros of some special functions. Also, the Laguerre-P\'olya class $\mathcal{LP}$ of real entire functions, which consist of uniform limits of real polynomials whose zeros are all real, was used intensively (for more details on the Laguerre-P\'olya class of entire functions we refer to \cite{BDM} and to the references therein). A real entire function $q$ belongs to the Laguerre-P\'{o}lya class $\mathcal{LP}$ if it can be represented in the form
$$
q(x) = c x^{m}e^{-\alpha x^{2}+\beta x} \prod_{n\geq1}
\left(1+\frac{x}{x_{n}}\right) e^{-\frac{x}{x_{n}}},
$$
where $c,$ $\beta,$ $x_{n}$ are real numbers, $\alpha \geq 0,$ $m$ is a natural number or zero, and $\sum\limits_{n\geq 1} x_{n}^{-2}$ converges. Motivated by the earlier works, in this work our main aim is to give some lower and upper bounds for the radii of convexity of some normalized $q$-Bessel functions. The results presented in this paper complement the results of \cite{BDM} about the radii of convexity and extend the known results from \cite{aktas2} on classical Bessel functions of the first kind to $q$-Bessel functions. In this study we consider two different normalized forms of Jackson and Hahn-Exton $q$-Bessel functions which are analytic in the unit disk of the complex plane. Because the functions $J_{\nu}^{(2)}(.;q)$ and $J_{\nu}^{(3)}(.;q)$ do not belong to $\mathcal{A}$, first we consider the following four normalized forms as in \cite{BDM}. For $\nu>-1$,
$$g_{\nu}^{(2)}(z;q)=2^{\nu}c_{\nu}(q)z^{1-\nu}J_{\nu}^{(2)}(z;q)=\sum_{n\geq0}\frac{(-1)^nq^{n(n+\nu)}}{4^n(q;q)_{n}(q^{\nu+1};q)_{n}}z^{2n+1},$$
$$h_{\nu}^{(2)}(z;q)=2^{\nu}c_{\nu}(q)z^{1-\frac{\nu}{2}}J_{\nu}^{(2)}(\sqrt{z};q)=\sum_{n\geq0}\frac{(-1)^nq^{n(n+\nu)}}{4^n(q;q)_{n}(q^{\nu+1};q)_{n}}z^{n+1},$$
$$g_{\nu}^{(3)}(z;q)=c_{\nu}(q)z^{1-\nu}J_{\nu}^{(3)}(z;q)=\sum_{n\geq0}\frac{(-1)^nq^{\frac{1}{2}n(n+1)}}{(q;q)_{n}(q^{\nu+1};q)_{n}}z^{2n+1},$$
$$h_{\nu}^{(3)}(z;q)=c_{\nu}(q)z^{1-\frac{\nu}{2}}J_{\nu}^{(3)}(\sqrt{z};q)=\sum_{n\geq0}\frac{(-1)^nq^{\frac{1}{2}n(n+1)}}{(q;q)_{n}(q^{\nu+1};q)_{n}}z^{n+1},$$ where $c_{\nu}(q)=(q;q)_{\infty}\big/(q^{\nu+1};q)_{\infty}$. As a result of the above normalizations, all of the above functions belong to the class $\mathcal{A}$.

\section{Bounds for the radii of convexity of some normalized $q$-Bessel functions}
In this section we give some tight lower and upper bounds for the radii of convexity of the above mentioned four normalized forms of the Jackson and Hahn-Exton $q$-Bessel functions. Also, we show that the radii of convexity of the above functions are solutions of some transcendental inequalities.

\begin{theorem}\label{th1}
Let $\nu>-1.$ Then the radius of convexity $r^{c}\left(g_{\nu}^{(2)}(z;q)\right)$ of the function $$z\mapsto g_{\nu}^{(2)}(z;q)=2^{\nu}c_{\nu}(q)z^{1-\nu}J_{\nu}^{(2)}(z;q)$$ is the smallest positive root of the equation $$(1-\nu)^2 J_{\nu}^{(2)}(r;q)+(3-2\nu)rdJ_{\nu}^{(2)}(r;q)/dr+r^2d^2J_{\nu}^{(2)}(r;q)/dr^2=0$$ and satisfies the following inequalities
$$\sqrt{\frac{4\left(1-q^{\nu+1}\right)\left(1-q\right)}{9q^{\nu+1}}} <r^{c}\left(g_{\nu}^{(2)}(z;q)\right)<\sqrt{\frac{36(q^2-1)\left(1-q^{\nu+1}\right)\left(1-q^{\nu+2}\right)}{q^{\nu+1}S_{\nu}(q)}},$$ $$2\sqrt[4]{\frac{(1+q)(1-q)^2\left(1-q^{\nu+1}\right)^2\left(q^{\nu+2}-1\right)}{q^{2(\nu+1)}S_{\nu}(q)}}<r^{c}\left(g_{\nu}^{(2)}(z;q)\right)<\sqrt{\frac{4\left(1-q^{\nu+1}\right)\left(1-q^{\nu+3}\right)T(q)S_{\nu}(q)}{q^{\nu+1}(1+q)\left(P_{\nu}(q)+R_{\nu}(q)\right)}},$$
where
$$P_{\nu}(q)=1458q-729q^{\nu+2}-1512q^{\nu+3}-2241q^{\nu+4}-837q^{\nu+5}-54q^{\nu+6}+479q^{\nu+7}+729q^{2\nu+5},$$
$$R_{\nu}(q)=783q^{2\nu+6}+783q^{2\nu+7}+152q^{2\nu+8}+98q^6-675q^4+54q^3+783q^2+729,$$
$$S_{\nu}(q)=31q^{\nu+3}+81q^{\nu+2}+50q^2-81q-81$$ and 
$$T(q)=(q-1)(q^3+2q^2+2q+1).$$
\end{theorem}
Note that by multiplying by $(1-q)^{-1}$ both sides of the above inequalities and taking the limit as $q\nearrow1$ for $\nu>-1$ we obtain the first two inequalities of \cite[Theorem 6]{aktas2}, namely:
\begin{equation}\label{eq2.1}
{\frac{2\sqrt{\nu+1}}{3}}<r^{c}(g_{\nu})<6\sqrt{\frac{(\nu+1)(\nu+2)}{56\nu+137}}
\end{equation}
and
\begin{equation}\label{eq2.2}
2\sqrt[4]{\frac{(\nu+1)^2(\nu+2)}{56\nu+137}}<r^{c}(g_{\nu})<\sqrt{\frac{2(\nu+1)(\nu+3)(56\nu+137)}{208\nu^2+1172\nu+1693}},
\end{equation}
where $r^{c}(g_{\nu})$ stands for the radii of convexity of the normalized Bessel function $$z\mapsto g_{\nu}(z)=2^{\nu}\Gamma(\nu+1)z^{1-\nu}J_{\nu}(z).$$

\begin{theorem}\label{th2}
Let $\nu>-1.$ Then the radius of convexity $r^{c}\left(h_{\nu}^{(2)}(z;q)\right)$ of the function $$z\mapsto h_{\nu}^{(2)}(z;q)=2^{\nu}c_{\nu}(q)z^{1-\frac{\nu}{2}}J_{\nu}^{(2)}(\sqrt{z};q)$$ is the smallest positive root of the equation 
$$\left(2-\nu\right)^2{J_{\nu}^{(2)}(\sqrt{r};q)}+\left(5-2\nu\right){\sqrt{r}}{dJ_{\nu}^{(2)}(\sqrt{r};q)/dr}+{r}{d^2J_{\nu}^{(2)}(\sqrt{r};q)/dr^2}=0$$ 
and satisfies the following inequalities $$\frac{(1-q)(1-q^{\nu+1})}{q^{\nu+1}}<r^c\left(h_{\nu}^{(2)}(z;q)\right)<
\frac{8(q^{\nu+1}-1)(q^{\nu+2}-1)(1-q^2)}{q^{\nu+1}U_{\nu}(q)},$$
$$\sqrt{\frac{8(1-q)^2(1+q)(1-q^{\nu+1})^2(1-q^{\nu+2})}{q^{2\nu+2}U_{\nu}(q)}}<
r^c\left(h_{\nu}^{(2)}(z;q)\right)<\frac{2(1-q^{\nu+1})(q^{\nu+3}-1)U_{\nu}(q)T(q)}{(1+q)q^{\nu+1}\left(M_{\nu}(q)+N_{\nu}(q)\right)},$$ where $$U_{\nu}(q)=\left(8q-8q^{\nu+2}+q^{\nu+3}-9q^2+8\right),$$
$$M_{\nu}(q)=32q-16q^{\nu+2}-21q^{\nu+3}-37q^{\nu+4}+6q^{\nu+5}+11q^{\nu+6}+3q^{\nu+7}$$ and
$$N_{\nu}(q)=16q^{2\nu+5}+5q^{2\nu+6}+5q^{2\nu+7}+q^{2\nu+8}+5q^2-11q^3-27q^4+12q^6+16.$$
\end{theorem}

Here we would like to emphasize that by multiplying by $(1-q)^{-2}$ both sides of the above inequalities and taking the limit as $q\nearrow1$ for $\nu>-1$ we obtain the first two inequalities of \cite[Theorem 7]{aktas2}, namely:
\begin{equation}\label{eq2.3}
\nu+1<r^{c}(h_{\nu})<\frac{16(\nu+1)(\nu+2)}{7\nu+23}
\end{equation}	
and
\begin{equation}\label{eq2.4}
\sqrt\frac{16(\nu+1)^2(\nu+2)}{{7\nu+23}}<r^{c}(h_{\nu})<\frac{2(\nu+1)(\nu+3)(7\nu+23)}{9\nu^2+60\nu+115},
\end{equation}
where $r^{c}(h_{\nu})$ stands for the radii of convexity of the normalized Bessel function  $$z\mapsto h_{\nu}(z)=2^{\nu}\Gamma(\nu+1)z^{1-\frac{\nu}{2}}J_{\nu}(\sqrt{z}).$$

\begin{theorem}\label{th3}
Let $\nu>-1.$ Then the radius of convexity $r^{c}\left(g_{\nu}^{(3)}(z;q)\right)$ of the function $$z\mapsto g_{\nu}^{(3)}(z;q)=c_{\nu}(q)z^{1-\nu}J_{\nu}^{(3)}(z;q)$$ is the smallest positive root of the equation $$(1-\nu)^2 J_{\nu}^{(3)}(r;q)+(3-2\nu) r dJ_{\nu}^{(3)}(r;q)/dr+r^2 d^2J_{\nu}^{(3)}(r;q)/dr^2=0$$ and satisfies the following inequalities
$$\sqrt{\frac{\left(1-q^{\nu+1}\right)\left(1-q\right)}{9q}} <r^{c}\left(g_{\nu}^{(3)}(z;q)\right)<\sqrt{\frac{9(q^2-1)\left(1-q^{\nu+1}\right)\left(1-q^{\nu+2}\right)}{qY_{\nu}(q)}},$$ $$\sqrt[4]{\frac{(1+q)(1-q)^2\left(1-q^{\nu+1}\right)^2\left(q^{\nu+2}-1\right)}{q^2Y_{\nu}(q)}}<r^{c}\left(g_{\nu}^{(3)}(z;q)\right)<\sqrt{\frac{\left(1-q^{\nu+1}\right)\left(1-q^{\nu+3}\right)T(q)Y_{\nu}(q)}{3q(q+1)\left(\theta_{\nu}(q)+\varphi_{\nu}(q)\right)}},$$
where
$$\theta_{\nu}(q)=261q-18q^{\nu+2}-504q^{\nu+3}-620q^{\nu+4}-504q^{\nu+5}-18q^{\nu+6}+67q^{2\nu+5},$$
$$\varphi_{\nu}(q)=261q^{2\nu+6}+261q^{2\nu+7}+243q^{2\nu+8}+67q^3+261q^2+243$$ and
$$Y_{\nu}(q)=81q^{\nu+3}+31q^{\nu+2}-31q-81.$$
\end{theorem}

Note that by multiplying by $(1-q)^{-1}$ both sides of the above inequalities and taking the limit as $q\nearrow1$ for $\nu>-1$ we obtain the following inequalities
\begin{equation}\label{eq2.5}
{\frac{\sqrt{\nu+1}}{3}}<r^{c}(g_{\nu})<3\sqrt{\frac{(\nu+1)(\nu+2)}{56\nu+137}}
\end{equation}
and
\begin{equation}\label{eq2.6}
\sqrt[4]{\frac{(\nu+1)^2(\nu+2)}{56\nu+137}}<r^{c}(g_{\nu})<\sqrt{\frac{(\nu+1)(\nu+3)(56\nu+137)}{2\left(208\nu^2+1172\nu+1693\right)}}.
\end{equation}

\begin{theorem}\label{th4}
Let $\nu>-1.$ Then the radius of convexity $r^{c}\left(h_{\nu}^{(3)}(z;q)\right)$ of the function $$z\mapsto h_{\nu}^{(3)}(z;q)=c_{\nu}(q)z^{1-\frac{\nu}{2}}J_{\nu}^{(3)}(\sqrt{z};q)$$ is the smallest positive root of the equation 
$$\left(2-\nu\right)^2{J_{\nu}^{(3)}(\sqrt{r};q)}+\left(5-2\nu\right){\sqrt{r}}{dJ_{\nu}^{(3)}(\sqrt{r};q)/dr}+{r}{d^2J_{\nu}^{(3)}(\sqrt{r};q)/dr^2}=0$$ 
and satisfies the following inequalities $$\frac{(1-q)(1-q^{\nu+1})}{4q}<r^c\left(h_{\nu}^{(3)}(z;q)\right)<
\frac{2(q^{\nu+1}-1)(q^{\nu+2}-1)(q^2-1)}{qE_{\nu}(q)},$$
$$\sqrt{\frac{(1-q)^2(1+q)(1-q^{\nu+1})^2(q^{\nu+2}-1)}{2q^{2}E_{\nu}(q)}}<
r^c\left(h_{\nu}^{(3)}(z;q)\right)<\frac{(1-q^{\nu+1})(q^{\nu+3}-1)E_{\nu}(q)T(q)}{2q(1+q)\left(K_{\nu}(q)+L_{\nu}(q)\right)},$$ where $$E_{\nu}(q)=\left(8q^{\nu+3}-q^{\nu+2}+q-8\right),$$
$$K_{\nu}(q)=5q+11q^{\nu+2}-21q^{\nu+3}-34q^{\nu+4}-21q^{\nu+5}+11q^{\nu+6}$$ and
$$L_{\nu}(q)=q^{2\nu+5}+5q^{2\nu+6}+5q^{2\nu+7}+16q^{2\nu+8}+5q^2+q^3+16.$$
\end{theorem}

Here we would like to emphasize that by multiplying by $(1-q)^{-2}$ both sides of the above inequalities and taking the limit as $q\nearrow1$ for $\nu>-1,$ we obtain the next two inequalities
\begin{equation}\label{eq2.7}
\frac{\nu+1}{4}<r^{c}(h_{\nu})<\frac{4(\nu+1)(\nu+2)}{7\nu+23}
\end{equation}	
and
\begin{equation}\label{eq2.8}
\sqrt\frac{(\nu+1)^2(\nu+2)}{{7\nu+23}}<r^{c}(h_{\nu})<\frac{(\nu+1)(\nu+3)(7\nu+23)}{2\left(9\nu^2+60\nu+115\right)}.
\end{equation}

It is important to mention that by making a comparison among of above obtained inequalities we have that the left-hand sides of \eqref{eq2.5} and \eqref{eq2.6} are weaker than the left hand sides of \eqref{eq2.1} and \eqref{eq2.2}, respectively. However, the right-hand sides of \eqref{eq2.5} and \eqref{eq2.6}  improve the right-hand sides of \eqref{eq2.1} and \eqref{eq2.2}, respectively. On the other hand, the left-hand sides of \eqref{eq2.7} and \eqref{eq2.8} are weaker than the left-hand sides of \eqref{eq2.3} and \eqref{eq2.4}, while the right-hand sides of \eqref{eq2.7} and \eqref{eq2.8} improve the right-hand sides of \eqref{eq2.3} and \eqref{eq2.4}.

\section{Proofs of main results}
\setcounter{equation}{0}

In this section we are going to present the proofs of our main results.

\begin{proof}[\bf Proof of Theorem \ref{th1}]
	
	By using the Alexander duality theorem for starlike and convex functions we can say that the function $g_{\nu}^{(2)}(z;q)$ is convex if and only if $z\mapsto{z\left(g_{\nu}^{(2)}(z;q)\right)^{\prime}}$ is starlike. But, the smallest positive zero of $z\mapsto{\left(z\left(g_{\nu}^{(2)}(z;q)\right)^{\prime}\right)^{\prime}}$ is actually the radius of starlikeness of $z\mapsto{z\left(g_{\nu}^{(2)}(z;q)\right)^{\prime}}$, according to \cite{bsk,bos}. Therefore, the radius of convexity $r^{c}(g_{\nu}^{(2)})$ is the smallest positive root of the equation $\left(r\left(g_{\nu}^{(2)}(r;q)\right)^{\prime}\right)^{\prime}=0.$ Thus, we get that the radius of convexity of the function $z\mapsto g_{\nu}^{(2)}(z;q)$ is the smallest positive root of the equation
	$$(1-\nu)^2 J_{\nu}^{(2)}(r;q)+(3-2\nu) r dJ_{\nu}^{(2)}(r;q)/dr+r^2 d^2J_{\nu}^{(2)}(r;q)/dr^2=0.$$
	
	Now, recall that the zeros $j_{\nu,n}(q), n\in\mathbb{N},$ of the Jackson $q$-Bessel function are all real and simple, according to \cite[Theorem 4.2]{ismail1}. Then, the function $g_{\nu}^{(2)}(\cdot;q)$ belongs to the Laguerre-P\'{o}lya class $\mathcal{LP}$ of real entire functions. Because of the properties of the class $\mathcal{LP}$ the function $z\mapsto {\left(z\left(g_{\nu}^{(2)}(z;q)\right)^{\prime}\right)^{\prime}}$ belongs also to the class $\mathcal{LP}.$ Hence the function $z\mapsto {\left(z\left(g_{\nu}^{(2)}(z;q)\right)^{\prime}\right)^{\prime}}$ has only real zeros. Also its growth order $\rho$ is $0$, that is $$\rho=\lim_{n\rightarrow\infty}\frac{n\log{n}}{2n\log2+\log(q;q)_{n}+\log(q^{\nu+1};q)_{n}-2\log(2n+1)-n(n+\nu)\log{q}}=0,$$
	since as $n\rightarrow\infty$ we have $(q;q)_{n}\rightarrow(q;q)_{\infty}<\infty \text{ and } (q^{\nu+1};q)_{n}\rightarrow(q^{\nu+1};q)_{\infty}<\infty$.
	Now, by applying Hadamard's Theorem \cite[p. 26]{Levin} we obtain
	$$G_{\nu}(z;q)={\left(z\left(g_{\nu}^{(2)}(z;q)\right)^{\prime}\right)^{\prime}}=\prod_{n\geq1}\left(1-\frac{z^2}{\left(\alpha_{\nu,n}(q)\right)^2}\right),$$ where ${\alpha_{\nu,n}(q)}$ is the $n$th zero of the function $G_\nu(\cdot;q)$. Now, via logarithmic derivation of $G_\nu(\cdot;q)$ we obtain
	\begin{equation}\label{eq3.1}
	\frac{G_\nu^{\prime}(z;q)}{G_\nu(z;q)}=-2\sum_{k\geq0}\epsilon_{k+1}z^{2k+1}, |z|<{\alpha_{\nu,1}(q)},
	\end{equation}
	where $\epsilon_k=\sum_{n\geq1}\left(\alpha_{\nu,n}(q)\right)^{-2k}$.
	Also, by using the infinite sum representation of $G_{\nu}$ we get
	\begin{equation}\label{eq3.2}
	\frac{G_\nu^{\prime}(z;q)}{G_\nu(z;q)}=\sum_{n\geq0}A_{n}z^{2n+1}\Bigg/\sum_{n\geq0}B_{n}z^{2n},
	\end{equation}
	where $$A_{n}=\frac{(-1)^{n+1}(2n+2)(2n+3)^2q^{(n+1)(n+\nu+1)}}{2^{2n+2}(q;q)_{n+1}(q^{\nu+1};q)_{n+1}}\ \ \ \mbox{and}\ \ \  B_{n}=\frac{(-1)^{n}(2n+1)^2q^{n(n+\nu)}}{2^{2n}(q;q)_{n}(q^{\nu+1};q)_{n}}.$$ By comparing \eqref{eq3.1} and \eqref{eq3.2} and matching all terms with the same degree we have the following Euler-Rayleigh sums $\epsilon_k=\sum_{n\geq1}\alpha_{\nu,n}^{-2k}(q)$ in terms of $\nu$ and $q$. That is, $$\epsilon_1=\frac{9q^{\nu+1}}{4(q^{\nu+1}-1)(q-1)},$$
	$$\epsilon_2=\frac{q^{2\nu+2}S_{\nu}(q)}{16(q^{\nu+1}-1)^2(q^{\nu+2}-1)(q-1)^2(q+1)}$$ and
	$$\epsilon_3=\frac{q^{3(\nu+1)}\left(P_{\nu}(q)+R_{\nu}(q)\right)}{64(q^{\nu+1}-1)^3(1-q^{\nu+2})(1-q^{\nu+3})(q-1)^2T(q)}.$$
	Now, by considering these Euler-Rayleigh sums in the known Euler-Rayleigh inequalities $$\epsilon_{k}^{-\frac{1}{k}}<\left(\alpha_{\nu,1}(q)\right)^2<\frac{\epsilon_{k}}{\epsilon_{k+1}}$$ for $\nu>-1$ and $k\in\{1,2\}$ we obtain the following inequalities 
	$$\sqrt{\frac{4\left(1-q^{\nu+1}\right)\left(1-q\right)}{9q^{\nu+1}}} <r^{c}\left(g_{\nu}^{(2)}(z;q)\right)<\sqrt{\frac{36(q^2-1)\left(1-q^{\nu+1}\right)\left(1-q^{\nu+2}\right)}{q^{\nu+1}S_{\nu}(q)}}$$ and $$2\sqrt[4]{\frac{(1+q)(1-q)^2\left(1-q^{\nu+1}\right)^2\left(q^{\nu+2}-1\right)}{q^{2(\nu+1)}S_{\nu}(q)}}<r^{c}\left(g_{\nu}^{(2)}(z;q)\right)<\sqrt{\frac{4\left(1-q^{\nu+1}\right)\left(1-q^{\nu+3}\right)T(q)S_{\nu}(q)}{q^{\nu+1}(1+q)\left(P_{\nu}(q)+R_{\nu}(q)\right)}}.$$
\end{proof}

\begin{proof}[\bf Proof of Theorem \ref{th2}]
	
	We proceed exactly as in the proof of Theorem \ref{th1}. By using the Alexander duality theorem for starlike and convex functions we can say that the function $h_{\nu}^{(2)}(z;q)$ is convex if and only if $z\mapsto{z\left(h_{\nu}^{(2)}(z;q)\right)^{\prime}}$ is starlike. But, the smallest positive zero of $z\mapsto{\left(z\left(h_{\nu}^{(2)}(z;q)\right)^{\prime}\right)^{\prime}}$ is actually the radius of starlikeness of $z\mapsto{z\left(h_{\nu}^{(2)}(z;q)\right)^{\prime}}$, according to \cite{bsk,bos}. Therefore, the radius of convexity $r^{c}(h_{\nu}^{(2)})$ is the smallest positive root of the equation $\left(r\left(h_{\nu}^{(2)}(r;q)\right)^{\prime}\right)^{\prime}=0.$ Thus, we get that the radius of convexity of the function $z\mapsto h_{\nu}^{(2)}(z;q)$ is the smallest positive root of the equation
	$$\left(2-\nu\right)^2{J_{\nu}^{(2)}(\sqrt{r};q)}+\left(5-2\nu\right){\sqrt{r}}{dJ_{\nu}^{(2)}(\sqrt{r};q)/dr}+{r}{d^2J_{\nu}^{(2)}(\sqrt{r};q)/dr^2}=0.$$
	
	Now, recall that the zeros $j_{\nu,n}(q), n\in\mathbb{N},$ of the Jackson $q$-Bessel function are all real and simple, according to \cite[Theorem 4.2]{ismail1}. Then, the function $h_{\nu}^{(2)}(\cdot;q)$ belongs to the Laguerre-P\'{o}lya class $\mathcal{LP}$ of real entire functions. Because of the properties of the class $\mathcal{LP}$ the function $z\mapsto {\left(z\left(h_{\nu}^{(2)}(z;q)\right)^{\prime}\right)^{\prime}}$ belongs also to the class $\mathcal{LP}.$ Hence the function $z\mapsto {\left(z\left(h_{\nu}^{(2)}(z;q)\right)^{\prime}\right)^{\prime}}$ has only real zeros. Also its growth order $\rho$ is $0$, that is $$\rho=\lim_{n\rightarrow\infty}\frac{n\log{n}}{2n\log2+\log(q;q)_{n}+\log(q^{\nu+1};q)_{n}-2\log(n+1)-n(n+\nu)\log{q}}=0,$$
	since as $n\rightarrow\infty$ we have $(q;q)_{n}\rightarrow(q;q)_{\infty}<\infty \text{ and } (q^{\nu+1};q)_{n}\rightarrow(q^{\nu+1};q)_{\infty}<\infty$.
	Now, by applying Hadamard's Theorem \cite[p. 26]{Levin} we obtain
	$$\Phi_{\nu}(z;q)={\left(z\left(h_{\nu}^{(2)}(z;q)\right)^{\prime}\right)^{\prime}}=\prod_{n\geq1}\left(1-\frac{z}{\beta_{\nu,n}(q)}\right),$$ where ${\beta_{\nu,n}(q)}$ is the $n$th zero of the function $\Phi_\nu(\cdot;q)$. Now, via logarithmic derivation of $\Phi_\nu(\cdot;q)$ we obtain
	\begin{equation}\label{eq3.3}
	\frac{\Phi_\nu^{\prime}(z;q)}{\Phi_\nu(z;q)}=-\sum_{k\geq0}\mu_{k+1}z^{k}, |z|<{\beta_{\nu,1}(q)},
	\end{equation}
	where $\mu_k=\sum_{n\geq1}\left(\beta_{\nu,n}(q)\right)^{-k}$.
	Also, by using the infinite sum representation of $\Phi_{\nu}$ we get
	\begin{equation}\label{eq3.4}
	\frac{\Phi_\nu^{\prime}(z;q)}{\Phi_\nu(z;q)}=\sum_{n\geq0}C_{n}z^{n}\Bigg/\sum_{n\geq0}D_{n}z^{n},
	\end{equation}
	where $$C_{n}=\frac{(-1)^{n+1}(n+1)(n+2)^2q^{(n+1)(n+\nu+1)}}{2^{2n+2}(q;q)_{n+1}(q^{\nu+1};q)_{n+1}}\ \ \ \mbox{and}\ \ \  D_{n}=\frac{(-1)^{n}(2n+1)^2q^{n(n+\nu)}}{2^{2n}(q;q)_{n}(q^{\nu+1};q)_{n}}.$$ By comparing \eqref{eq3.3} and \eqref{eq3.4} and matching all terms with the same degree we have the following Euler-Rayleigh sums $\mu_k=\sum_{n\geq1}\beta_{\nu,n}^{-k}(q)$ in terms of $\nu$ and $q$. That is, $$\mu_1=\frac{q^{\nu+1}}{(q^{\nu+1}-1)(q-1)},$$
	$$\mu_2=\frac{q^{2\nu+2}U_{\nu}(q)}{8(q^{\nu+1}-1)^2(1-q^{\nu+2})(q-1)^2(q+1)}$$ and
	$$\mu_3=\frac{q^{3(\nu+1)}\left(M_{\nu}(q)+N_{\nu}(q)\right)}{16(q^{\nu+1}-1)^3(q^{\nu+2}-1)(q^{\nu+3}-1)(q-1)^2T(q)}.$$
	Now, by considering these Euler-Rayleigh sums in the known Euler-Rayleigh inequalities $$\mu_{k}^{-\frac{1}{k}}<\beta_{\nu,1}(q)<\frac{\mu_{k}}{\mu_{k+1}}$$ for $\nu>-1$ and $k\in\{1,2\}$ we obtain the next inequalities
	$$\frac{(1-q)(1-q^{\nu+1})}{q^{\nu+1}}<r^c\left(h_{\nu}^{(2)}(z;q)\right)<
	\frac{8(q^{\nu+1}-1)(q^{\nu+2}-1)(1-q^2)}{q^{\nu+1}U_{\nu}(q)},$$
	$$\sqrt{\frac{8(1-q)^2(1+q)(1-q^{\nu+1})^2(1-q^{\nu+2})}{q^{2\nu+2}U_{\nu}(q)}}<
	r^c\left(h_{\nu}^{(2)}(z;q)\right)<\frac{2(1-q^{\nu+1})(q^{\nu+3}-1)U_{\nu}(q)T(q)}{(1+q)q^{\nu+1}\left(M_{\nu}(q)+N_{\nu}(q)\right)}.$$
\end{proof}

\begin{proof}[\bf Proof of Theorem \ref{th3}]
	
	By using the Alexander duality theorem for starlike and convex functions we can say that the function $g_{\nu}^{(3)}(z;q)$ is convex if and only if $z\mapsto{z\left(g_{\nu}^{(3)}(z;q)\right)^{\prime}}$ is starlike. But, the smallest positive zero of $z\mapsto{\left(z\left(g_{\nu}^{(3)}(z;q)\right)^{\prime}\right)^{\prime}}$ is actually the radius of starlikeness of $z\mapsto{z\left(g_{\nu}^{(3)}(z;q)\right)^{\prime}}$, according to \cite{bsk,bos}. Therefore, the radius of convexity $r^{c}(g_{\nu}^{(3)})$ is the smallest positive root of the equation $\left(r\left(g_{\nu}^{(3)}(r;q)\right)^{\prime}\right)^{\prime}=0.$ Thus, we get that the radius of convexity of the function $z\mapsto g_{\nu}^{(3)}(z;q)$ is the smallest positive root of the equation
	$$(1-\nu)^2 J_{\nu}^{(3)}(r;q)+(3-2\nu)rdJ_{\nu}^{(3)}(r;q)/dr+r^2 d^2J_{\nu}^{(3)}(r;q)/dr^2=0.$$
	
	Now, recall that the zeros $l_{\nu,n}(q), n\in\mathbb{N},$ of the Hahn-Exton $q$-Bessel function are all real and simple, according to \cite[Theorem 4.2]{ismail1}. Then, the function $g_{\nu}^{(3)}(\cdot;q)$ belongs to the Laguerre-P\'{o}lya class $\mathcal{LP}$ of real entire functions. Because of the properties of the class $\mathcal{LP}$ the function $z\mapsto {\left(z\left(g_{\nu}^{(3)}(z;q)\right)^{\prime}\right)^{\prime}}$ belongs also to the class $\mathcal{LP}.$ Hence the function $z\mapsto {\left(z\left(g_{\nu}^{(3)}(z;q)\right)^{\prime}\right)^{\prime}}$ has only real zeros. Also its growth order $\rho$ is $0$, that is $$\rho=\lim_{n\rightarrow\infty}\frac{n\log{n}}{\log(q;q)_{n}+\log(q^{\nu+1};q)_{n}-2\log(2n+1)-\frac{n(n+1)}{2}\log{q}}=0,$$
	since as $n\rightarrow\infty$ we have $(q;q)_{n}\rightarrow(q;q)_{\infty}<\infty \text{ and } (q^{\nu+1};q)_{n}\rightarrow(q^{\nu+1};q)_{\infty}<\infty$.
	Now, by applying Hadamard's Theorem \cite[p. 26]{Levin} we obtain
	$$\mathcal{H}_{\nu}(z;q)={\left(z\left(g_{\nu}^{(3)}(z;q)\right)^{\prime}\right)^{\prime}}=\prod_{n\geq1}\left(1-\frac{z^2}{\left(h_{\nu,n}(q)\right)^2}\right),$$ where ${h_{\nu,n}(q)}$ is the $n$th zero of the function $\mathcal{H}_\nu(\cdot;q)$. Now, via logarithmic derivation of $\mathcal{H}_\nu(\cdot;q)$ we obtain
	\begin{equation}\label{eq3.5}
	\frac{\mathcal{H}_\nu^{\prime}(z;q)}{\mathcal{H}_\nu(z;q)}=-2\sum_{k\geq0}\eta_{k+1}z^{2k+1}, |z|<{h_{\nu,1}(q)},
	\end{equation}
	where $\eta_k=\sum_{n\geq1}\left(h_{\nu,n}(q)\right)^{-2k}$.
	Also, by using the infinite sum representation of $\mathcal{H}_{\nu}$ we have
	\begin{equation}\label{eq3.6}
	\frac{\mathcal{H}_\nu^{\prime}(z;q)}{\mathcal{H}_\nu(z;q)}=\sum_{n\geq0}E_{n}z^{2n+1}\Bigg/\sum_{n\geq0}F_{n}z^{2n},
	\end{equation}
	where $$E_{n}=\frac{(-1)^{n+1}(2n+2)(2n+3)^2q^{\frac{(n+1)(n+2)}{2}}}{(q;q)_{n+1}(q^{\nu+1};q)_{n+1}}\ \ \ \mbox{and}\ \ \  F_{n}=\frac{(-1)^{n}(2n+1)^2q^{\frac{n(n+1)}{2}}}{(q;q)_{n}(q^{\nu+1};q)_{n}}.$$ By comparing \eqref{eq3.5} and \eqref{eq3.6} and matching all terms with the same degree we have the following Euler-Rayleigh sums $\eta_k=\sum_{n\geq1}h_{\nu,n}^{-2k}(q)$ in terms of $\nu$ and $q$. That is, $$\eta_1=\frac{9q}{(q^{\nu+1}-1)(q-1)},$$
	$$\eta_2=\frac{q^{2}Y_{\nu}(q)}{(q^{\nu+1}-1)^2(q^{\nu+2}-1)(q-1)^2(q+1)}$$ and
	$$\eta_3=\frac{3q^{3}\left(\theta_{\nu}(q)+\varphi_{\nu}(q)\right)}{64(q^{\nu+1}-1)^3(1-q^{\nu+2})(1-q^{\nu+3})(q-1)^2T(q)}.$$
	Now, by considering these Euler-Rayleigh sums in the known Euler-Rayleigh inequalities $$\eta_{k}^{-\frac{1}{k}}<\left(h_{\nu,1}(q)\right)^2<\frac{\eta_{k}}{\eta_{k+1}}$$ for $\nu>-1$ and $k\in\{1,2\}$ we obtain the following inequalities
	$$\sqrt{\frac{\left(1-q^{\nu+1}\right)\left(1-q\right)}{9q}} <r^{c}\left(g_{\nu}^{(3)}(z;q)\right)<\sqrt{\frac{9(q^2-1)\left(1-q^{\nu+1}\right)\left(1-q^{\nu+2}\right)}{qY_{\nu}(q)}},$$ $$\sqrt[4]{\frac{(1+q)(1-q)^2\left(1-q^{\nu+1}\right)^2\left(q^{\nu+2}-1\right)}{q^2Y_{\nu}(q)}}<r^{c}\left(g_{\nu}^{(3)}(z;q)\right)<\sqrt{\frac{\left(1-q^{\nu+1}\right)\left(1-q^{\nu+3}\right)T(q)Y_{\nu}(q)}{3q(q+1)\left(\theta_{\nu}(q)+\varphi_{\nu}(q)\right)}}.$$
\end{proof}

\begin{proof}[\bf Proof of Theorem \ref{th4}]
	
	By using the Alexander duality theorem for starlike and convex functions we can say that the function $h_{\nu}^{(3)}(z;q)$ is convex if and only if $z\mapsto{z\left(h_{\nu}^{(3)}(z;q)\right)^{\prime}}$ is starlike. But, the smallest positive zero of $z\mapsto{\left(z\left(h_{\nu}^{(3)}(z;q)\right)^{\prime}\right)^{\prime}}$ is actually the radius of starlikeness of $z\mapsto{z\left(h_{\nu}^{(3)}(z;q)\right)^{\prime}}$, according to \cite{bsk,bos}. Therefore, the radius of convexity $r^{c}(h_{\nu}^{(3)})$ is the smallest positive root of the equation $\left(r\left(h_{\nu}^{(3)}(r;q)\right)^{\prime}\right)^{\prime}=0.$ Thus, we get that the radius of convexity of the function $z\mapsto h_{\nu}^{(3)}(z;q)$ is the smallest positive root of the equation
$$\left(2-\nu\right)^2{J_{\nu}^{(3)}(\sqrt{r};q)}+\left(5-2\nu\right){\sqrt{r}}{dJ_{\nu}^{(3)}(\sqrt{r};q)/dr}+{r}{d^2J_{\nu}^{(3)}(\sqrt{r};q)/dr^2}=0.$$
	
	Now, recall that the zeros $l_{\nu,n}(q), n\in\mathbb{N},$ of the Hahn-Exton $q$-Bessel function are all real and simple, according to \cite[Theorem 4.2]{ismail1}. Then, the function $h_{\nu}^{(3)}(\cdot;q)$ belongs to the Laguerre-P\'{o}lya class $\mathcal{LP}$ of real entire functions. Because of the properties of the class $\mathcal{LP}$ the function $z\mapsto {\left(z\left(h_{\nu}^{(3)}(z;q)\right)^{\prime}\right)^{\prime}}$ belongs also to the class $\mathcal{LP}.$ Hence, the function $z\mapsto {\left(z\left(h_{\nu}^{(3)}(z;q)\right)^{\prime}\right)^{\prime}}$ has only real zeros. Also its growth order $\rho$ is $0$, that is $$\rho=\lim_{n\rightarrow\infty}\frac{n\log{n}}{\log(q;q)_{n}+\log(q^{\nu+1};q)_{n}-2\log(n+1)-\frac{n(n+1)}{2}\log{q}}=0,$$
	since as $n\rightarrow\infty$ we have $(q;q)_{n}\rightarrow(q;q)_{\infty}<\infty \text{ and } (q^{\nu+1};q)_{n}\rightarrow(q^{\nu+1};q)_{\infty}<\infty$.
	Now, by applying Hadamard's Theorem \cite[p. 26]{Levin} we obtain
	$$\psi_{\nu}(z;q)={\left(z\left(h_{\nu}^{(3)}(z;q)\right)^{\prime}\right)^{\prime}}=\prod_{n\geq1}\left(1-\frac{z}{\gamma_{\nu,n}(q)}\right),$$ where ${\gamma_{\nu,n}(q)}$ is the $n$th zero of the function $\psi_\nu(\cdot;q)$. Now, via logarithmic derivation of $\psi_\nu(\cdot;q)$ we obtain
	\begin{equation}\label{eq3.7}
	\frac{\psi_\nu^{\prime}(z;q)}{\psi_\nu(z;q)}=-\sum_{k\geq0}\sigma_{k+1}z^{k}, |z|<{\gamma_{\nu,1}(q)},
	\end{equation}
	where $\sigma_k=\sum_{n\geq1}\left(\gamma_{\nu,n}(q)\right)^{-k}$.
	Also, by using the infinite sum representation of $\Phi_{\nu}$ we get
	\begin{equation}\label{eq3.8}
	\frac{\psi_\nu^{\prime}(z;q)}{\psi_\nu(z;q)}=\sum_{n\geq0}G_{n}z^{n}\Bigg/\sum_{n\geq0}H_{n}z^{n},
	\end{equation}
	where $$G_{n}=\frac{(-1)^{n+1}(n+1)(n+2)^2q^{\frac{(n+1)(n+2)}{2}}}{(q;q)_{n+1}(q^{\nu+1};q)_{n+1}}\ \ \ \mbox{and}\ \ \  H_{n}=\frac{(-1)^{n}(n+1)^2q^{\frac{n(n+1)}{2}}}{(q;q)_{n}(q^{\nu+1};q)_{n}}.$$ By comparing \eqref{eq3.7} and \eqref{eq3.8} and matching all terms with the same degree we have the following Euler-Rayleigh sums $\sigma_k=\sum_{n\geq1}\gamma_{\nu,n}^{-k}(q)$ in terms of $\nu$ and $q$. That is, $$\sigma_1=\frac{4q}{(q^{\nu+1}-1)(q-1)},$$
	$$\sigma_2=\frac{2q^{2}E_{\nu}(q)}{(q^{\nu+1}-1)^2(1-q^{\nu+2})(1-q)^2(1-q)}$$ and
	$$\sigma_3=\frac{4q^{3}\left(K_{\nu}(q)+L_{\nu}(q)\right)}{(1-q^{\nu+1})^3(1-q^{\nu+2})(1-q^{\nu+3})(1-q)^2T(q)}.$$
	Now, by considering these Euler-Rayleigh sums in the known Euler-Rayleigh inequalities $$\sigma_{k}^{-\frac{1}{k}}<\gamma_{\nu,1}(q)<\frac{\sigma_{k}}{\sigma_{k+1}}$$ for $\nu>-1$ and $k\in\{1,2\}$ we obtain the next two inequalities
	$$\frac{(1-q)(1-q^{\nu+1})}{4q}<r^c\left(h_{\nu}^{(3)}(z;q)\right)<
	\frac{2(q^{\nu+1}-1)(q^{\nu+2}-1)(q^2-1)}{qE_{\nu}(q)},$$
	$$\sqrt{\frac{(1-q)^2(1+q)(1-q^{\nu+1})^2(q^{\nu+2}-1)}{2q^{2}E_{\nu}(q)}}<
	r^c\left(h_{\nu}^{(3)}(z;q)\right)<\frac{(1-q^{\nu+1})(q^{\nu+3}-1)E_{\nu}(q)T(q)}{2q(1+q)\left(K_{\nu}(q)+L_{\nu}(q)\right)}.$$
\end{proof}

\end{document}